\documentclass[12pt]{amsart}
\def\bitem#1#2#3#4{\bibitem{#1} {\sc #2}, {\it #3}, {\rm #4}.}

\def\C{\mathbb C}

\let\secc\section

\def\Re{\operatorname{\mathrm{Re}}}
\def\Im{\operatorname{\mathrm{Im}}}
\let\colon:

\newcounter{thrm}

\newtheorem{corollary}[thrm]{Corollary}
\newtheorem{definition}[thrm]{Definition}

\newtheorem{theorem}[thrm]{Theorem}

\theoremstyle{definition}
\newtheorem{remark}[thrm]{Remark}

\def\mybib{\section*{References}\small\list
{\arabic{enumi}.}{\settowidth\labelwidth{[99]}\leftmargin\labelwidth
\advance\leftmargin\labelsep \usecounter{enumi}} }

\title{Reducing Essential Eigenvalues in the Boundary of the Numerical
Range{ }}

\author{Norberto Salinas}
\address{
\hskip-\parindent
Norberto Salinas\\
Department of Mathematics\\
University of Kansas\\
Lawrence KS 66045}
\email{norberto@@kuhub.cc.ukans.edu}

\author{Maria Victoria Velasco}
\address{
\hskip-\parindent
Maria Victoria Velasco\\
Department of Analysis\\
University of Granada\\
Granada, Spain}
\email{vvelasco@@goliat.ugr.es}

\thanks{Research at MSRI is supported in part by MSRI grant
DMS-9022140.}

\begin{document}
\begin{abstract}
We examine a purely geometric property of a point in the boundary of
the numerical range of a (Hilbert space) operator that implies that
such a point is a reducing essential eigenvalue of the given operator.
Roughly speaking, such a property means that the boundary curve of the
numerical range has infinite curvature at that point (we must exclude
however linear verteces because they may be reducing eigenvalues
without being reducing essential eigenvalues). This result allows us
to give an elegant proof of a conjecture of Joel Anderson: {\it A
compact perturbation of a scalar multiple of the identity operator can
not have the closure of its numerical range equal to half a disk
(neither equal to any acute circular sector).}
\end{abstract}

\maketitle

In this note, we examine a purely geometric
property of a point
in the boundary of the numerical range of a (Hilbert space) operator 
which implies that such a point is a reducing essential eigenvalue of
the given
operator (see Theorem \ref{princi1} and Corollary \ref{princi2}). Roughly speaking, such a
property means
that the boundary curve of the numerical range has infinite
curvature at that
point (we must exclude however linear verteces  because they may be
reducing
eigenvalues without being reducing essential eigenvalues). This
result
allows us to give an elegant proof of a conjecture of Joel Anderson
(see \cite{Ag}, \cite{Lan}, \cite{Ra}):

 \begin{quote} {\it   A compact
perturbation of a scalar
multiple of the identity operator can not have the closure of its
numerical range
equal to half a disk (neither equal to any acute circular sector).}
\end{quote}

 Indeed, there are two non-smooth points in the boundary of
such a convex
set that have infinite curvature and that are not linear verteces.
So they are reducing essential eigenvalues (from Corollary 8)
which is not
possible for a multiple of the identity operator.

  The
above mentioned theorem also proves in the affirmative a
strengthening of a
conjecture of Mathias H\"ubner \cite{Hu}:  such points (of
infinite
curvature) are in the essential spectrum of the given operator.

\secc{Reducing Essential Eigenvalues and  the Numerical Range}

In this section, we recall some basic facts about the numerical range of an operator on Hilbert space, and the notion of reducing essential eigenvalues.  
\begin{definition} {\rm Let $T$
be an operator
acting on a (fixed)
Hilbert space ${\mathcal H}$. The {\it numerical range} and the
{\it state numerical range} of $T$ are defined, respectively, as follows
\[ W(T) := \{ (\langle T x,
x\rangle
\colon \; x\in {\mathcal H}, \;
\|x\| =1 \}, \] \[ W_s(T) :=\{ \phi
T\colon  \; \phi \in  S(C^*(T)) \},\] where $C^*(T)$
denotes the $C^*$-algebra generated by $T$ and the identity
operator $I$ on
${\mathcal H}$.  Also, we used the notation $S({\mathcal A})$ for the
state-space of the
$C^*$-algebra ${\mathcal A}$.} \end{definition}

 \begin{remark}  
\begin{enumerate}
\item[(a)]
We recall \cite{St-Wi} that $W(T)$ is always
convex, and $W_s(T) = \overline{W(T)}$.
\item[(b)]
The problem of determining the class of  bounded convex sets of the form $W(T)$ for some operator $T$ on ${\mathcal H}$ is still open (an easy cardinality argument shows that there are many bounded convex sets which are not in that class). In \cite{Ag} it is shown that $W(T)$ is always a borel subset of
 ${\C }$, and in \cite{Ra} it is proved that the components  of $W(T)$
 in the boundary of $W(T)$ must be  either singletons or conic arcs.  On the
other hand, it is not even known whether $W(T)$ can be the union of
the open unit disk and a non-trivial  open arc in its boundary. 
\end{enumerate}
\end{remark}  

 \begin{definition} {\rm Let
 $\lambda $ be a point in the
boundary of a convex subset $C$ of the complex plane  (with non-empty
 interior). We say that $\lambda$ is of {\it infinite
curvature} if, after a suitable rotation and translation which identifies
$\lambda $ with zero, we can choose the real axis as supporting line for
$C$  (in which case, we shall say that $C$ is in standard position) and then
$$
\lim_{\stackrel{\Re \alpha \to 0}{\stackrel{ \Im \alpha \to 0}{\alpha \in
C}}} \frac{\vert \Im
\alpha \vert}{ \Re^2 \alpha} =  \infty .
$$
For convenience sake, we adopt the convention that
 $C$ is in standard position if it is contained in the upper-half plane,
and the point of the boundary under study is the origin.
If the real axis is the only tangent line to the boundary of $C$, then
 we call the origin a point of {\it infinite lefthand curvature} whenever
$$
\lim_{\stackrel{\Re \alpha \to 0-}{\stackrel{ \Im \alpha \to 0}{\alpha \in
C}}} \frac{\vert \Im
\alpha \vert}{ \Re^2 \alpha} = \infty ,
$$ and we call the origin a point of {\it infinite righthand curvature}
 whenever
$$
\lim_{\stackrel{\Re \alpha \to 0+}{\stackrel{ \Im \alpha \to 0}{\alpha \in
C}}} \frac{\vert \Im
\alpha \vert}{ \Re^2 \alpha} = \infty .
$$}
\end{definition}

\begin{remark} 
\indent (a)
Obviously, in order to show that a point is of infinite curvature, the
suitable rotation will not be unique in general. Of course, this will
not be an obstacle for deducing later some analitic properties that such
 a points may satisfy.
Moreover,  it should be noticed that  the property
of having infinite curvature depends only on the graph of
the
function defined by the boundary of $C$ near $\lambda =0$. Indeed,
let $C$ be in a standard position. Then it is
easily shown that $\lambda$ is of infinite
curvature precisely when
$$
\lim_{\stackrel{ \alpha \to 0}{ \alpha \in 
\partial C}} \frac{ \Im
\alpha }{ \Re^2 \alpha} = \infty .
$$
\indent (b) Of course, a
linear vertex in the boundary of a convex set $C$  (i.e., a point $\lambda$ in 
$\partial C$ such that $C$ can be written as the convex hull of $\lambda$ and
 another convex set $C'\not = C$)  is a point of
infinite curvature of $C$.
 In fact, it readily follows that even a corner in the boundary of
a convex set
is a point of infinite curvature. (We recall \cite{Lan} that a point $\lambda$ in
the boundary of a convex set $C$ is called a {\it corner} if $C$ is a
subset of a convex set
$C'$ where $\lambda$ is a linear vertex). But, it is easy to come up
with examples
of points of infinite curvature where the boundary is continuously
differentiable. \newline
\indent(c) To have  infinite lefthand curvature means that
$C$, and more precisely the graph given by $
\partial C$, has infinite curvature on the left of zero, while to
 have  infinite righthand curvature means that
the curvature on the right  of zero is infinite.
It is also clear that a point may be of infinite
 onesided curvature but not of infinite curvature. For instance, this
 is the case for the origin in the convex set  $C := A \cup B$ where
$A := \{\alpha  \in {\C } \colon \ -1 \leq
\Re \alpha \leq 0, \ 0 \leq \Im \alpha \leq 1 \}$ and $B:=   \{ \alpha  \in  {\C}\colon
\ ( \Re \alpha)^{3/2} \leq  \Im \alpha \leq 1, \ 0 \leq \Re \alpha \leq 1 \}$.
\end{remark}

 \begin{definition} {\rm A point $\lambda $ is called an
{\it reducing approximate  eigenvalue} of an operator $T$ if there exists
a sequence
$u_n$ of unit vectors in ${\mathcal H}$ such that \[
\lim_{n\rightarrow \infty } \|(T-\lambda I)
u_n\| + \|(T^*-\overline{\lambda} I) u_n\| =0. \] If in addition $u_n$
tends weakly
to zero, then $\lambda $ is called a {\it reducing essential
eigenvalue}. The set of all
reducing approximate eigenvalues of the operator $T$ will be
denoted by
$R(T)$ and the set of reducing essential eigenvalues by
$R_e(T)$. \newline
\indent A point
$\lambda \in {\C }$ is called a {\it normal essential eigenvalue} of
$T$ if $\lambda$ is in  $R_e(T)$
and every time that an orthonormal sequence $u_n$ satisfies
$\|(T- \lambda I)
u_n \| \to 0$, we also have $\|(T^*-\overline{\lambda }I)
u_n\|\to 0$}.  \end{definition}

 \begin{remark} \indent(a) The
properties and characterizations of $R(T)$ and $R_e(T)$ were
already
discussed in \cite{Sa}. In particular, $\lambda$ belongs to $R(T)$ if, and only
if, there exist a $*$-homomorphism $\phi: C^*(T)\rightarrow {\C }$ such that $\phi
T=\lambda$
and, analogously, $\lambda$ is in $R_e(T)$ if and only if there exist a $*$-homomorphism
$\phi: C^*(T) \rightarrow {\C }$ such that $\phi(C^*(T) \cap {\mathcal K})=0$
and  $\phi T=\lambda$.
 (Here, we are employing the usual notation
${\mathcal K}$ for the ideal of compat operators)
Thus, $R(T)$ (resp. $R_e(T)$)
is contained in the intersection
of $\sigma (T)$  (resp. $\sigma_{e} (T)$) and the complex conjugate set of
$\sigma(T^*)$ (resp. $\sigma_e (T^*)$).
 \newline
\indent (b) Let $R_{00}(T)$ be
the set of finite
multiplicity reducing eigenvalues of $T$ which are isolated points
of $R(T)$. Here, by a finite-dimensional reducing eigenvalue of $T$ we
mean a complex number $\lambda$ such that $null (T - \lambda I) \cap null
(T^* - \overline{\lambda}I)$ is a non trivial finite dimensional subspace.
Then, it is shown in \cite{Sa}, Theorem 6.1, that
$R(T) = R_e (T) \cup  R_{00}(T)$, where the union (of course) is disjoint.
\end{remark}

 \begin{theorem}\label{princi1} \begin{enumerate}
\item[(a)] Let $ 0 \in 
\partial W_s(T)$ be a point of infinite
curvature and consider $W_s(T)$ in standard position. If $u_n$ is a
sequence of unit vectors in ${\mathcal H}$ such that  $\langle T u_n,
u_n\rangle \rightarrow 0$  then $\| Tu_n\| \to 0$ and $\| T^*u_n\| \to 0$.
\newline
\item[(b)] Let $ 0 \in 
\partial W_s(T)$ be a point of infinite righthand (resp. lefthand)
curvature and consider $W_s(T)$ in standard position. If $u_n$ is a
sequence of unit vectors in ${\mathcal H}$ such that  $\langle T u_n,
u_n\rangle \rightarrow 0$ where, for all n,
$\langle T u_n,
u_n\rangle$ is not zero and is contained in a segment of end points
 given by zero and
some $\alpha_0$ in $W_s(T)$ with $\Re \alpha_0 >0$ (resp. $\Re \alpha_0
<0$), then $\| Tu_n\| \to 0$ and $\| T^*u_n\| \to 0$. 
\end{enumerate}
 \end{theorem} 

 \begin{proof}
First of all we will stablish some inequalities which will be useful to our
purposes.
We begin by
observing that from the condition $W_s(T) \subset \{\lambda \in {\C}\colon \ \Im \lambda \geq 0 \}$ we
obtain the following elementary but important inequality: \[ \Im
\langle Tu, u\rangle
\geq 0, \ \mbox{ for all } u \mbox{ in } {\mathcal H}. \] Also, given
arbitrary vectors
$u,z$ in ${\mathcal H}$ and an arbitrary $\alpha\in  {\C }$, we have:
\small
\begin{eqnarray*} \langle T( u + \alpha z),
u + \alpha z\rangle  =  \langle Tu,u\rangle  +
\overline{\alpha} \langle Tu, z\rangle  + \alpha \langle Tz,
u\rangle   + \vert  \alpha \vert^2 \langle Tz,z\rangle .
\end{eqnarray*}
\normalsize
 Since $\Im \langle T( u \pm \alpha z), u
\pm \alpha z\rangle  \geq 0$, we see that
\small
 \begin{eqnarray}
\vert \Im [ \overline{\alpha} \langle Tu, z\rangle
+ \alpha \langle Tz, u\rangle  ] \vert\leq \nonumber \\ \Im [
\langle Tu,u\rangle  + \vert  \alpha
\vert^2 \langle Tz,z\rangle  ] \leq \vert \langle Tu,u\rangle
\vert +  \vert  \alpha \vert^2 \|z\|^2 \| T\|.
\end{eqnarray}
\normalsize
Consequently,
\begin{eqnarray}
\vert \Im \langle T (u + \alpha z),
u + \alpha z\rangle  \vert \leq 2
 \Im [ \langle Tu,u\rangle  + \vert  \alpha \vert^2 \langle
Tz,z\rangle  ] \leq  \nonumber \\  2 (\vert
\langle Tu,u\rangle  \vert + \vert  \alpha \vert^2 \| T\| \|z\|^2 ).
\end{eqnarray} On the other hand, \begin{eqnarray}
\vert \Re \langle T( u + \alpha z), u + \alpha z\rangle  \vert \geq
\nonumber \\ \vert
\Re [ \overline{\alpha} \langle Tu, z\rangle  + \alpha \langle Tz,
u\rangle  ] \vert -\vert \Re \langle Tu,u\rangle
\vert - \vert \alpha \vert^2 \vert \Re \langle Tz,z\rangle
\vert \geq \nonumber \\ \vert \Re [ \overline{\alpha} \langle Tu,
z\rangle  + \alpha
\langle Tz, u\rangle  ] \vert -
\vert \langle Tu,u\rangle  \vert -  \vert \alpha \vert^2 \vert
\|T\| \|z\|^2.
\end{eqnarray}

Now we treat the case where zero is a point of infinite curvature
in $
\partial W_s(T)$, {\it i.e.},
assertion (a). Let $u_n$ be a sequence of unit vectors
in ${\mathcal H}$ such that $\langle  Tu_n, u_n\rangle\rightarrow 0$.
We may write
(uniquely) $Tu_n =\delta _n u_n +
\beta _n v_n$ and  $T^* u_n = \overline{\delta _n } u_n + \overline{\gamma }_n w_n$, where
$v_n$ and $w_n$ are unit vectors
orthogonal to $u_n$, and $\beta _n, \gamma _n, \delta _n$, are given
respectively by $\delta_n : = \langle Tu_n, u_n\rangle , \ \beta_n : = \langle Tu_n,
v_n\rangle $ and $\gamma_n: = \langle Tw_n, u_n\rangle $, for all $n=1,2,
\ldots $. Also, by multiplying by an appropriate complex number
of modulus one, if needed, we may assume that $\Re \langle v_n, w_n\rangle  \geq 0$
for all $n$. Let $\tau_n$ be a complex number with $ \vert \tau_n \vert$=1 and such that
$\vert \overline{\tau_n }  \beta_n+ \tau_n \gamma_n\vert =\vert \beta_n \vert +
\vert \gamma_n \vert$, $n=1, 2, \ldots$. Let $r_n:= \vert \delta_n \vert^{\frac{1}{2}}$ if $\delta_n \neq 0$
and $r_n:= \frac{1}{n}$ otherwise. We define $\alpha_n:= r_n \tau_n$, for all n. Since
$\vert \overline{\alpha_n }  \beta_n
+ \alpha_n \gamma_n\vert = r_n (\vert \beta_n \vert + \vert \gamma_n\vert)$,
 there exist suitable 
 complex numbers $\eta_n$ with $\vert \eta_n\vert=1$
such that $\overline{\alpha_n }  \beta_n
+ \alpha_n \gamma_n =r_n (\vert \beta_n \vert + \vert \gamma_n \vert )
\eta_n$, for all $n$.
Let $z_n :=v_n + w_n$. Since
$\langle Tu_n, z_n\rangle =\beta_n (1+ \langle v_n, w_n\rangle )$ and
$\langle Tz_n,
u_n\rangle  = \gamma_n (1+ \langle v_n,
w_n\rangle )$ we obtain:
\begin{eqnarray*} \overline{\alpha_n }
\langle Tu_n, z_n\rangle  + \alpha_n
\langle Tz_n,u_n\rangle = ( \overline{\alpha_n }  \beta_n
+ \alpha_n \gamma_n) (1+ \langle v_n, w_n\rangle )= \\
 r_n
 (\vert \beta_n \vert + \vert \gamma_n\vert)\eta_n ( 1 +\langle v_n, w_n\rangle )
\end{eqnarray*}
Hence, from (1.1) we deduce that
$$
 r_n
 (\vert \beta_n \vert + \vert\gamma_n\vert)
\vert \Im [ \eta_n (1+ \langle v_n,w_n\rangle ) ] \vert \leq \vert
\delta_n
\vert + 4 {r_n}^2 \|T\|
$$
so, we obtain:
\begin{eqnarray}
(\vert \beta_n \vert + \vert \gamma_n\vert) \vert \Im [ \eta _n (1+
\langle v_n, w_n\rangle ) ] \vert \rightarrow 0.
\end{eqnarray}
We claim that also \begin{eqnarray}
(\vert \beta_n \vert + \vert \gamma_n\vert) \vert \Re [
\eta _n (1+
\langle v_n, w_n\rangle ) ] \vert \rightarrow 0.
\end{eqnarray}
To prove this claim
we define   $x_n := u_n + \alpha_n z_n$ whenever
$\Im [\overline{\alpha_n} \langle
Tu_n,  z_n \rangle +\alpha_n \langle Tz_n,u_n \rangle] \geq 0$ and
$x_n := u_n - \alpha_n z_n$ otherwise.
Now we consider two particular cases: \newline
\indent (i) In the case that $\Re
\langle Tx_n,x_n\rangle = 0$ for all $n$ it follows from (1.3) that
$r_n(\vert \beta_n \vert + \vert \gamma_n\vert)\vert \Re [\eta_n (1
+ \langle v_n, w_n \rangle ) ]\vert  \leq
\vert \delta_n \vert + 4{r_n}^2 \| T\|$. Thus $(\vert \beta_n \vert + \vert \gamma_n\vert)\vert \Re [\eta_n (1
+ \langle v_n, w_n \rangle ) ]\vert \to 0$, as desired.\newline
\indent (ii)  Assume now that $\Re \langle Tx_n,x_n\rangle \neq 0$ for all $n$.
Also  $\Im \langle Tx_n,x_n\rangle \neq 0$ because
the curvature
at $\lambda =0$ is infinite so that the only real point of $W(T)$ is zero.  Furthermore, it follows that
 \[ \mu _n: = \frac {\vert \Im \langle Tx_n, x_n\rangle  \vert}{ \Re^2\langle
Tx_n , x_n \rangle  }
\rightarrow \infty, \] since $ \langle Tx_n , x_n\rangle \to 0$ with
$\|x_n\| \to 1$ .
Therefore by $(1.2)$
\begin{eqnarray*} \vert \Re \langle Tx_n , x_n\rangle  \vert =
 \left( \frac{ \vert \Im \langle Tx_n , x_n \rangle  \vert }{ \mu _n }
\right)^{\frac{1}{2}} \leq    \nonumber \\
\left( \frac{2 ( \vert \delta_n \vert + 4  {r_n}^2
\|T\|)}{\mu _n } \right)^{\frac{1}{2}}.
\end{eqnarray*} However by $(1.3)$
\begin{eqnarray*} \vert \Re \langle Tx_n , x_n\rangle\vert
\geq r_n(\vert \beta_n \vert + \vert \gamma_n\vert)
 \vert \Re  \eta _n  (1+
\langle v_n , w_n \rangle ) ] \vert - \vert \delta_n
\vert - 4
{r_n}^2   \|T\|.
\end{eqnarray*} Hence,
\small
 \begin{eqnarray*}
r_n(\vert \beta_n \vert + \vert \gamma_n\vert)\vert  \Re [ \eta _n  (1+ \langle v_n ,
w_n \rangle ) ] \vert  \leq  \\
\left( \frac{ 2(\vert \delta_n \vert + 4 {r_n}^2\|T\|)}{\mu _n } \right)^{\frac{1}{2}} +
\vert \delta_n \vert+ 4 {r_n}^2 \| T\|
 \end{eqnarray*} and so $(\vert \beta_n \vert + \vert \gamma_n\vert)\vert \Re [\eta_n (1
+ \langle v_n, w_n \rangle ) ]\vert \to 0$, as desired. \newline
\normalsize
\noindent Thus, from (i) and (ii), the above claim follows. Therefore, by
(1.4), we see
 that
\begin{eqnarray*} (\vert \beta_n \vert + \vert \gamma_n\vert)\vert   (1+ \langle
v_n , w_n \rangle ) \vert = (\vert \beta_n \vert + \vert \gamma_n\vert)\vert
\eta_n (1+ \langle
v_n , w_n \rangle ) \vert =
\\
(\vert \beta_n \vert + \vert \gamma_n\vert) \left(
\Re^2 [ \eta _n  (1+ \langle v_n , w_n \rangle ) ]
 +
\Im^2 [ \eta _n  (1+ \langle v_n , w_n \rangle ) ]
 \right)^{\frac{1}{2}}
\rightarrow 0. \end{eqnarray*}
 But $\vert 1 + \langle v_n,w_n \rangle  \vert \geq 1$
because
$\Re \langle v_n, w_n\rangle \geq 0$, for every $n$. Thus $\vert \beta_n
\vert + \vert \gamma_n\vert \to 0$ which means
precisely that
$\| Tu_n  \| \to 0$ as well as $\| T^* u_n  \| \to 0$, so that (a)
is proved.   \newline
\indent In order to prove (b), assume for example that
zero is a point of
infinite righthand curvature (the case of infinite lefthand curvature may
be handled similarly). Let $\alpha_0$ be in $W_s(T)$ such that $\Re
\alpha_0 >0$. Also  $\Im \alpha_0 >0$ since otherwise zero would not
be a point of infinite righthand curvature. Let $\varepsilon_n$ be a
sequence of positive real numbers such that  $\varepsilon_n \to 0$  and let
$u_n$ be a sequence of unit vectors in ${\mathcal H}$ such that
$\langle Tu_n , u_n\rangle = \varepsilon_n \alpha_0$. As in the proof of
 (a), we write  $Tu_n =\delta _n u_n +
\beta _n v_n$ and  $T^* u_n = \overline{\delta _n } u_n + \overline{\gamma }_n w_n$, where
$v_n$ and $w_n$ are as before. However, in this case we have the following
additional information $\delta_n = \varepsilon_n \alpha_0$. Let $z_n := v_n
+ w_n$ and $M$ be a constant greater than
$\frac{ \vert \alpha_0 \vert \vert \Re \langle Tz_n,
z_n\rangle \vert}{\Re \alpha_0}$ (take for
instance
$M= \frac{5 \vert \alpha_0\vert \|T\| }{\Re \alpha_0}$). We define again
$\alpha_n : = r_n \tau_n$, for every $n$, where $\tau_n$ is as before but
now $r_n:= \left( \frac{\vert \delta_n \vert}{M}\right)^{\frac{1}{2}}$.
Then,
$$
\vert \alpha_n\vert^2 \vert \langle T z_n, z_n \rangle \vert =  \varepsilon_n \frac{ \vert \alpha_0 \vert  \vert \Re \langle Tz_n ,
z_n \rangle \vert}{M} < \varepsilon_n \Re \alpha_0 = \Re \langle Tu_n ,
u_n \rangle
$$
so that
$$\Re \langle Tu_n , u_n \rangle + \vert \alpha_n \vert^2 \Re \langle Tz_n ,
z_n \rangle >0.$$ On the other hand, for every n, there exists a 
complex number $\eta_n$ with $\vert \eta_n\vert =1$ such that
 \begin{eqnarray*} \overline{\alpha_n }
\langle Tu_n, z_n\rangle  + \alpha_n
\langle Tz_n,u_n\rangle =  \left( \frac{ \vert \delta_n \vert}{M}\right)^{\frac{1}{2}}
(\vert \beta_n \vert + \vert \gamma_n \vert ) \eta _n (1+ \langle v_n,
w_n\rangle ).
\end{eqnarray*}
Now we define $x_n:= u_n + \alpha_n z_n$ if
$\Re \eta _n (1+ \langle v_n, w_n\rangle ) >0$ and $x_n:= u_n - \alpha_n
z_n$ otherwise. Then, 
$\Re \langle Tx_n ,x_n \rangle >0$.
Consequently, since there are not positive real numbers in $W_s(T)$ we
deduce that
$\Im \langle Tx_n ,x_n \rangle \neq 0$ for all $n$.
Also, again because zero is a point of infinite righthand curvature and $\langle Tx_n
,x_n \rangle \to 0$ with $\| x_n\| \to 1$, it follows that
$\frac {\vert \Im \langle Tx_n, x_n\rangle  \vert}{ \Re^2\langle
Tx_n , x_n \rangle  }
\rightarrow \infty$. Then, by proceding as above in the corresponding part of the
proof of (a), we deduce that  $(\vert \beta_n \vert + \vert \gamma_n \vert )\vert \Re [\eta_n (1
+ \langle v_n, w_n \rangle ) ]\vert \to 0$.
On the other hand, from (1.2) we obtain $(\vert \beta_n \vert + \vert \gamma_n \vert )\vert \Im [\eta_n (1
+ \langle v_n, w_n \rangle ) ]\vert \to 0$. Therefore  we conclude that $(\vert \beta_n \vert + \vert \gamma_n \vert )\vert \eta_n (1
+ \langle v_n, w_n \rangle ) \vert \to 0$ and so  $\vert \beta_n \vert + \vert \gamma_n \vert
  \to 0$. This proves that  $\| Tu_n\| \to 0$ and $\| T^*u_n\| \to 0$.
\end{proof}

\indent Assertion (a) in the next corollary
constitutes an improvement of the main result in \cite{Hu}. Also
assertion (b) shows in particular  that {\it points of infinite one-sided
curvature of $
\partial W_s(T)$ which are not of infinite curvature are reducing essential
eigenvalues}. Finally the assertion (c)
improves \cite{Lan} , Corollary 4.


 \begin{corollary}\label{princi2}
\begin{enumerate}
\item[(a)] If $\lambda $ is either a point of infinite
curvature or a point of infinite one-sided curvature of $
\partial
W_s(T)$, then $\lambda \in R(T)$.  \newline \noindent
\item[(b)] If $\lambda$ is a
point of infinite one-sided curvature in $
\partial W_s(T)$ but it is
not a linear vertex, then $\lambda\in R_e(T)$.\newline
\item[(c)]
If $\lambda $ is a point of infinite curvature in $
\partial W_s(T)$, but
it is not a linear
vertex of $W(T)$, then  $\lambda $ is a
normal essential
eigenvalue of $T$.  
\end{enumerate}\end{corollary} 

 \begin{proof}
Part (a) is an immediate consequence of the
the definition of $R(T)$ and Theorem \ref{princi1}. In order to prove
 (b) and (c) we first claim that a
point $\lambda$ in $W(T)$ which is not a linear vertex belongs to $R_e(T)$.
To show the claim we  consider $T$ decomposed  in the (direct) sum
$T= A \oplus  B$, acting on $H= H_0^\bot \oplus H_0$  where
$H_0: = \mbox{Ker}(T - \lambda I) \cap \mbox{Ker}(T^* - \overline{\lambda}
I)$.
Since $A$ and $B$ are the restriction of $T$ to $H_0^\bot$ and $H_0$
respectively, it is clear that $\lambda$ is not a reducing eigenvalue of
$A$ so, in particular $\lambda\notin R_{00}(A)$. Notice that it is easily proved
that $W(T)= \mbox{Convhull}(\{\lambda\},W(A))$  (see also \cite{Ra}).
Hence $\lambda$ must be either in
$
\partial  W(A)$ or in the complement of ${W_s(A)}$.
In the first case, $W(T) = W(A)$. The
second alternative
is not possible since, in that case, $\lambda$ would be a linear vertex
of $W_s(T)$
which contradicts our working assumption. Since $\lambda$ is a point of
infinite
curvature in $
\partial  W_s(A)$ we deduce, from (a), that
$\lambda\in  R(A)$. Because $\lambda\notin  R_{00}(A)$, we conclude that
$\lambda\in R_e(A)$. Consequently $\lambda \in R_e(T)$ and the claim is
shown. To complete the proof of the theorem we merely observe that if
$u_n$ is a orthonormal sequence such that
$\| (T- \lambda I) u_n\| \to 0$ then $\langle Tu_n, u_n\rangle\to
\lambda$.  Hence, as a consequence of Theorem \ref{princi1}, we see that
also  $\|(T-\lambda)^*u_n\|\to 0$.
 \end{proof}

 The
following result proves a strengthening of a conjecture of Mathias
H\"ubner \cite{Hu}.

 \begin{corollary} Assume that $\lambda $ is a
point of infinite one-sided curvature
of $
\partial  W_s(T)$, for some operator $T$, and suppose that there
is only one
tangent line to $
\partial  W_s(T)$ at the point $\lambda $. Then
$\lambda \in  R_e(T)$. If, in
addition, $\lambda $ is a point of infinite curvature of $W_s(T)$
then $\lambda $ is a
normal essential eigenvalue of $T$. \end{corollary} 

 \begin{proof}  It follows immediately from
part (c) of the previous corollary
because $\lambda $ can not be a
linear vertex of $
\partial W_s(T)$.
 \end{proof}

\secc{Joint Reducing  Essential  Eigenvalues and 
the
Joint Numerical Range}

Now, we shall find a mild extension of the last corollary
to $n$-tuples.

 \begin{definition} {\rm Let
${\mathbf T}= (T_1, T_2, \ldots , T_n)$ be an $n$-tuple of operators
acting on
${\mathcal H}$. The {\it joint numerical range} and the {\it joint state numerical 
range}  
of ${\mathbf T}$ are defined, respectively, as follows:
\[ W({\mathbf T}) := \{ (\langle T_1 x,
x\rangle , \langle T_2 x, x\rangle ,
\ldots , \langle T_n x, x\rangle )\colon  \; x\in {\mathcal H}, \;
\|x\| =1 \}, \] \[ W_s({\mathbf T}) :=\{ (\phi
T_1, \phi  T_2, \ldots , \phi  T_n)\colon  \; \phi \in  S(C^*({\mathbf T})) \},\] where $C^*({\mathbf T})$
denotes the $C^*$-algebra generated by ${\mathbf T}$ and the identity
operator ${\mathbf I}$ on
${\mathcal H}$.
 }\end{definition}

 \begin{remark} \indent (a)
 $W_s({\mathbf T})$ is always
convex but, in general,  $W_s({\mathbf T}) = \overline{W({\mathbf T})}$ only for $n=1$.
 Indeed, for $n>1$ the above equality does not 
always hold and
$W({\mathbf T})$ is not convex in general \cite{Da}.
\newline
\indent (b) If  $\lambda $ is
an extreme point of $
\partial  W_s({\mathbf T})$, for a given $n$-tuple
of operators ${\mathbf T}$,
then $\lambda \in  \overline{W({\mathbf T})}$.
Indeed, the proof of this fact uses
a standard argument: let $\Sigma _\lambda $ be the subset of
$S(C^*({\mathbf T}))$ consisting of
those states $\phi $ on $C^*({\mathbf T})$ such that $\phi ({\mathbf T})=\lambda $. Since $\Sigma _\lambda $ is
compact and convex in the $w^*$-topology, there exists an extreme
point $\psi $ in
$\Sigma _\lambda $. Using the fact that $\lambda $ is an extreme
point of $W({\mathbf T})$, it
readily follows that $\psi $ is a pure state of $C^*({\mathbf T})$. By
Glim's Lemma (see \cite{Dix} and  \cite{Gli}), we obtain that $\psi $ is in
the $w^*$-closure
of the set of vector states of $C^*({\mathbf T})$. Thus, evaluating at
${\mathbf T}$, we
conclude that $\lambda \in  \overline{W({\mathbf T})}$, as desired.
\newline
\indent (c) In the next theorem we
use the notions of joint reducing approximate point spectrum and of
joint reducing essential spectrum of an $n$-tuple of operators.
These are the natural
extension to $n$-tuples of the corresponding notions for single
operators. \end{remark}

\begin{theorem} Given an $n$-tuple ${\mathbf T}=(T_1, \cdots T_n)$ of operators acting on
${\mathcal H}$,  assume that $\lambda=(\lambda_1, \cdots , \lambda_n) \in 
\partial  W_a({\mathbf T})$ satisfies
$\lambda _j$ is a point of
one-sided infinite curvature of $
\partial W_s(T_j)$, for $1\leq  j\leq  n$. Then
$\lambda \in  R({\mathbf T})$. If, furthermore, $\lambda $ is
not a linear vertex of $
\partial  W({\mathbf T})$, then $\lambda \in R_e({\mathbf T})$. If,
in addition, $\lambda _j$ is a point of infinite curvature of $
\partial
W_s(T_j)$, for $1\leq  j\leq  n$, then $\lambda $ is a joint normal essential eigenvalue for ${\mathbf T}$.
\end{theorem}

 \begin{proof} {\rm The proof
consits of a repeated application of Corollary \ref{princi2}, by
projecting in each
coordinate. We should point out that under the present hypotheses,
$\lambda $ is an
extreme point of $\overline{W({\mathbf T})}$}.
\end{proof}

 \begin{remark} \indent (a) It would be interesting to
weaken the hypotheses of Corollary \ref{princi2}. For example, is it possible to
simply assume
that $\lambda $ is a point of upper-infinite curvature in $
\partial
W_s(T)$ (in the sense that we replace $\lim$ by $\lim\sup$ in the definition of infinite curvature)? \newline
\indent (b) We recall that the $m$-th matricial range of an
operator ${\mathbf T}$ on
${\mathcal H}$ is the set of $m\times m$ matrices of the form $\phi
({\mathbf T})$, where $\phi \colon \ C^*({\mathbf T})\to {\mathcal M}_m$ is a unital completely positive linear map
(see, for instance, \cite{Bu-Sa}). Would it be possible to define a notion of
points of
infinite matricial curvature in the boundary of the matricial
range of ${\mathbf T}$,
so that to be able to conclude that such points are actually
in the
approximate reducing matricial spectrum of ${\mathbf T}$ (those $\phi ({\mathbf T})$
for which $\phi $
is actually a *-homomorphism)? \end{remark}

\begin{mybib}
\bitem{Ag}
{J. Agler}
{Geometric and topological properties of the
numerical range}
{Indiana Univ. Math. J. 31 (1982), 767-777}

\bitem{Bu-Sa}
{ J. Bunce $\&$ N. Salinas}
{Completely positive maps on
$C^*$-algebras and the left matricial spectra of an operator}
{Duke Math. J. 74 (1976), 747-773}

\bitem{Da}
{ A. T. Dash} { Joint numerical ranges}{Glasnik Matematicki
(1972) 75-81}

\bitem{Dix}
{ J. Dixmier}{  Les $C^*$-alg\`{e}bres at leurs repres\'{e}ntations}
{Gauthier-Villars, Paris (1964)}

\bitem{Gli}
{ J. Glimm}{  A Stone-Weierstrass theorem for $C^*$-algebras}
{Ann. Math. 72 (1960), 216-244}

\bitem{Hu}
{ M. H$\ddot{u}$bner}{ Spectrum where the boundary of the
numerical range is not round}{ To appear in Rocky Mountain Journal of
Math.}

\bitem{Lan}
{ J. Lancaster}{  The boundary of the numerical range}
{Proc. Amer. Math. Soc. 49 (1975), 393-398}

\bitem{Ra}
{ M. Radjabalipour $\&$ H. Radjavi}{ On the geometry of numerical
ranges}{ Pacific J. Math.  61 (1975), 507-511}

\bitem{Sa}
{ N. Salinas}{ Reducing essential eigenvalues}
{Duke. Math. J. 40 (1973), 561-580}

\bitem{St-Wi}
{ J. G. Stampfli $\&$ J.P. Willians}{ Growth conditions and the
numerical range in a Banach algebra}{ T$\hat{o}$hoku Math. Jour. 20 (1968),
417-424}

\end{mybib}

\end{document}